\documentclass[12pt]{amsart}
\usepackage{amssymb}
\usepackage{amsmath,latexsym}
\usepackage{hyperref} 
\usepackage{eucal}
\newcommand{\rr}{\mathbb{R}}
\newcommand{\p}{\partial}

\newcommand{\cc}{\mathbb{C}}
\newtheorem{corollary}{Corollary}[section]
\newtheorem{lemma}{Lemma}[section]

\newtheorem{theorem}{Theorem}[section]

\newtheorem*{TA}{Theorem A}
\newtheorem*{TB}{Theorem B}
\newtheorem*{TC}{Theorem C}
\newtheorem*{TD}{Theorem D}

\title[KdV Equations]
{On decay  properties of solutions of 
the $k$-generalized KdV equations}
\author{Pedro Isaza}
\address[P. Isaza]{Departamento  de Matem\'aticas\\
Universidad Nacional de Colombia\\ A. A. 3840, Medellin\\Colombia}
\email{pisaza@unal.edu.co}
\author{Felipe Linares}
\address[F. Linares]{IMPA\\
Instituto Matem\'atica Pura e Aplicada\\
Estrada Dona Castorina 110\\
22460-320, Rio de Janeiro, RJ\\Brazil}
\email{linares@impa.br}

\author{Gustavo Ponce}
\address[G. Ponce]{Department  of Mathematics\\
University of California\\
Santa Barbara, CA 93106\\
USA.}
\email{ponce@math.ucsb.edu}
\keywords{Korteweg-de Vries  equation,  weighted Sobolev spaces}
\subjclass{Primary: 35Q53. Secondary: 35B05}

\begin{document}
\begin{abstract} We prove  special decay properties of solutions to the initial value problem associated to the $k$-generalized  Korteweg-de Vries  equation.
These are  related with  persistence properties of the solution flow in weighted Sobolev spaces and with   sharp unique
continuation properties of solutions to this equation. As application of our method we also obtain results concerning the decay behavior of perturbations of the traveling wave solutions as well as results for solutions corresponding to special data.

\end{abstract}


\maketitle

\numberwithin{equation}{section}

\section{Introduction}

In this work we shall study special decay properties of real solutions to the initial value problem (IVP)
associated to the  $k$-generalized Korteweg-de Vries ($k$-gKdV) equation 
\begin{equation}
\label{kgKdV}
\begin{cases}
\partial_t u +\partial_x^3u +u^k\partial_x u = 0,\;\;\;\;\;\; x\in\rr,\;\;k\in\mathbb Z^+,\\
u(x,0) = u_0(x).
\end{cases}
\end{equation}

 These decay properties of the solution $u(t)$ will be measured in appropriate weighted 
 $L^2(w  dx)$-spaces. 
 
 First, we shall be concerned with asymmetric (increasing) weights, for which the result
 will be restricted to forward times $t>0$. In this regard we find the result in \cite{Ka} for the KdV equation, $k=1$ in \eqref{kgKdV}, in the space $L^2(e^{\beta x}dx),\,\beta>0$. There it was shown that the persistence property holds for $L^2$-solutions in $L^2(e^{\beta x}dx),\,\beta>0$, for $t>0$ (persistence property in the function space $X$ means that the solution $u(\cdot)$ describes a continuous curve on $X$, $u\in C([0,T]:X)$). Moreover, formally in this space the operator $\partial_t+\partial_x^3$  becomes $\partial_t+(\partial_x-\beta)^3$ so the solutions of the equation exhibit a parabolic behavior. More precisely, the following result for the KdV equation was proven in \cite{Ka} (Theorem 11.1 and Theorem 12.1)

   \begin{TA} \label{Ka}
  Let $ u\in C([0,\infty)\,:\,H^2(\rr))$ be a solution of the IVP \eqref{kgKdV} with $k=1$ and 
  \begin{equation}
  \label{1.3}
  u_0\in H^2(\rr)\cap L^2(e^{\beta x}dx),\;\;\;\;\text{for some}\;\;\beta>0,
 \end{equation}
  then
  \begin{equation}
  \label{1.4} 
e^{\beta x}u\in   C([0,\infty)\,:\,L^2(\rr))\cap C((0,\infty)\,:\,H^{\infty}(\rr)),
\end{equation}
with
\begin{equation}
\label{1.5}
\|u(t)\|_2=\|u_0\|_2,\;\;\;\;\;\;\;\|u(t)-u_0\|_{-3,2}\leq Kt,\;\;\;\;\;\;t>0,
\end{equation}
\begin{equation}
\label{1.6}
\|e^{\beta x}u(t)\|_2\leq e^{Kt}\,\| e^{\beta x}u_0\|_2,\;\;\;\;\;t>0,
\end{equation}
  and
 \begin{equation}
\label{1.7}
\int_0^{\infty} \,e^{-Kt}\,\| e^{\beta x}\partial_xu(t)\|_2^2dt \leq \frac{1}{4\beta} \|e^{\beta x}u_0\|_2^2,
\end{equation}
where $K=K(\beta,\|u_0\|_2)$. 

Moreover, the map data-solution $u_0\to u(t)$ is continuous from $L^2(\rr)\cap  L^2(e^{\beta x}dx)$ to  $C([0,T]\,:\,L^2(e^{\beta x}dx))$, for any $T>0$.
 \end{TA}

On a similar regard, in \cite{EKPV06} a unique continuation result was established. This gives an upper bound for the possible space decay of solutions to the IVP \eqref{kgKdV}:

\begin{TB} \label{EKPV06}
 There exists $c_0>0$ such that for any pair  
 $$
u_1,\,u_2\in C([0,1]:H^4(\rr)\cap L^2(|x|^2dx))
$$
of  solutions of  \eqref{kgKdV}, if 
\begin{equation}
 \label{3:2}
 u_1(\cdot,0)-u_2(\cdot,0),\,\;\, u_1(\cdot,1)-u_2(\cdot,1)\in L^2(e^{c_0x_{+}^{3/2}}dx),
\end{equation}  then $u_1\equiv u_2$.
  \end{TB}
 
 Above we used the notation:   $ x_{+}=\max\{x;\,0\}$. Similarly, we will use later on  $ x_{-}=\max\{-x;\,0\}$
 

The power $3/2$ in the exponent in \eqref{3:2} reflects the asymptotic behavior of the Airy function. 
The solution 
of the initial value problem (IVP) 
\begin{equation}
 \begin{aligned}
 \begin{cases}
 \partial_t v + \partial_x^3 v=0,\\
 v(x,0)=v_0(x),
 \end{cases}
 \end{aligned}
\end{equation}
is given by the group $\{U(t)\,:\,t\in \rr\}$
 $$
 U(t)v_0(x)=\frac{1}{\root{3}\of{3t}}\,Ai\left(\frac{\cdot}{\root {3}\of{3t}}\right)\ast v_0(x), 
 $$
where
 $$
Ai(x)=c\,\int_{-\infty}^{\infty}\,e^{ ix\xi+i \xi^3/3 }\,d\xi,
 $$
 is the Airy function which satisfies the estimate 
$$
 |Ai(x)|\leq c\,\frac{e^{-c x_{+}^{3/2}}}{(1+x_{-})^{1/4}}.
 $$

Observe that Theorem B gives a restriction on the possible decay 
of a non-trivial solution of 
\eqref{kgKdV}
at two different times. More precisely, 
taking $u_2\equiv 0$ one has that if $u_1(t)$ is a solution of the IVP \eqref{kgKdV} such that
\begin{equation}\label{a1}
|u_1(x,t)|\leq e^{-a_0x_{+}^{3/2}}\;\;\text{at}\;\;t=0,\,1,\;\;\text{for}\;\;a_0>>1, \;\text{then\hskip10pt }u_1\equiv 0.
\end{equation}

Our first theorem shows that the above result is close to be optimal.  By rescaling, the persistence property can not hold in the space $L^2(e^{a_0 x_{+}^{3/2}}dx)$ in an arbitrary large time interval. However,  it does it with a factor $a(t)$ in front of the exponential term $x_{+}^{3/2}$ which measures how this exponential property decreases with time. 

To simplify the exposition, we shall first state our result  in the case of the KdV, i.e.  $k=1$ in \eqref{kgKdV} :

 \begin{theorem} \label{A1}
Let $a_0$ be a positive constant. For any given data
\begin{equation}
\label{0.1}
u_0\in L^2(\rr)\cap L^2(e^{a_0x_{+}^{3/2}}dx),
\end{equation}
the unique solution of the IVP \eqref{kgKdV} provided by Theorem C below satisfies that for any $T>0$
\begin{equation}
\label{0.2}
\sup_{t\in [0,T]}\,\int_{-\infty}^{\infty} e^{a(t)x_{+}^{3/2}}|u(x,t)|^2dx  \leq C^*=C^*( a_0, \|u_0\|_2, \|e^{a_0x_{+}^{3/2}/2}u_0\|_2, T),
\end{equation}
with
\begin{equation}
 \label{0.3}
a(t)=\frac{a_0}{(1+27 a_0^2t/4)^{1/2}}.
 \end{equation}

\end{theorem}

This result extends to the difference of two appropriate solutions of the IVP \eqref{kgKdV} with $k=1$.

 \begin{theorem} \label{B1}
Let $a_0$ be a positive constant. Let $u(t),\,v(t)$ be solutions of the IVP \eqref{kgKdV} with $k=1$ such that
\begin{equation}
\label{4.1}
\begin{aligned}
& u\in C([0,\infty)\,:\,H^1(\rr)\cap L^2(|x|dx)) \cap\dots, \\
& v\in C([0,\infty)\,:\,H^1(\rr)) \cap\dots 
\end{aligned}
\end{equation}
If
\begin{equation}
\label{4.2}
\int_{-\infty}^{\infty} e^{a_0x_{+}^{3/2}}|u(x,0)-v(x,0)|^2dx  <\infty,
\end{equation}
then for any $T>0$
\begin{equation}
\label{4.3}
\sup_{t\in[0,T]}\,\int_{-\infty}^{\infty} e^{a(t)x_{+}^{3/2}}|u(x,t)-v(x,t)|^2dx  \leq C^*,
\end{equation}
with 
\begin{equation}
\label{4.41}
C^*=C^*( a_0, \|u_0\|_{1,2}, \|v_0\|_{1,2}, \| |x|^{1/2}u_0\|_2, \|u_0-v_0\|_{1,2}, \|e^{a_0x_{+}^{3/2}/2}(u_0-v_0)\|_2, T),
\end{equation}
 and 
\begin{equation}
 \label{4.42}
a(t)=\frac{a_0}{(1+27 a_0^2t/4)^{1/2}}.
 \end{equation}

\end{theorem}

The results in Theorem \ref{A1} and Theorem \ref{B1} apply to other powers $k$ in the IVP \eqref{kgKdV}:
\begin{theorem} \label{BB1}
Let $a_0$ be a positive constant. For any given data 
$$
u_0\in H^{s_k}(\rr)\cap L^2(e^{a_0x_{+}^{3/2}})
$$
with
$\,s_2=1/4$ if $k=2$, $s_3=0$ if $k=3$, and $s_k>(k-4)/2k,\;k\geq 4$, then the unique solution of the IVP \eqref{kgKdV} provided 
by Theorem  D below satisfies \eqref{0.2} with $a(t)$ as in \eqref{0.3} for any $T>0$ if $k=2,3$ and for $T=T(\|u_0\|_{s_k,2})>0$ for $k\geq 4$.
\end{theorem}

Similarly, for Theorem \ref{B1}.

Remarks:

(a) Following the argument in the proof of Theorem 1.4 in \cite{EKPV06} one has  that given $a_0>0,\,\epsilon>0$, there exist  $u_0\in \mathbb S(\rr)$, $c_1,\,c_2>0$ and $\Delta T>0$ such that  the corresponding solution $u(x,t)$ of the IVP \eqref{kgKdV} with $k=1$
satisfies
$$
c_1\,e^{-(a_0+\epsilon)x^{3/2}}\leq u(x,t)\leq c_2\,e^{-(a_0-\epsilon)x^{3/2}},\;\;\;\;\;\,x>>1,\;\;t\in [0,\Delta T].
$$

(b) For different powers $k$ in \eqref{kgKdV} one has that the same statement of Theorem \ref{A1} is valid except for
the last part concerning the continuity of the map data-solution $u_0\to u(t)$ which will be continuous
from $H^{s_k}(\rr)\cap L^2(e^{\beta x}dx)$ to $C([0,T]: L^2(e^{\beta x}dx))$ with $s_k$ and $T$ as in the statement of Theorem \ref{BB1}.

(c)  For other results involving asymmetric weights of a polynomial type we refer to \cite{KrFa} and \cite{GiTs}.

 (d) Consider the 1-D semi-linear Schr\"odinger equation
 \begin{equation}
\label{NLS}
\partial_t v =i( \partial_x^2 v + F(v,\overline v)),
\end{equation}
with $F:\cc^2\to \cc$, $F\in C^{2}$  and $F(0)=\partial_uF(0)=\partial_{\bar
u}F(0)=0$.
As far as we are aware  it is unknown whether or not there exist non-trivial solutions $v(t)$ of \eqref{NLS} satisfying
$$
 v(t)\in L^2(e^{a_0x_+^{1+\epsilon}}\,)\;\;\;\;\;\,t\in[0,T],\;\;\;\;\;\;\;\text{for some}\;\;T>0,\;\,\epsilon>0.
 $$

Next, we consider weighted spaces with symmetric weight of the form
$$
L^2(\langle x\rangle^bdx)=L^2((1+x^2)^{b/2}dx)
$$
 for which persistent properties should hold regardless of the time direction considered, i.e. forward $t>0$ or backward $t<0$.

In this setting our second result establishes that for a solution of the IVP \eqref{kgKdV} to satisfy the persistent property in  $L^2(\langle x\rangle^bdx)$ it needs to have a similar property in an appropriate Sobolev space $H^s(\rr)$, i.e. decay in $L^2$ can only hold if $u(t)$ is regular enough in $L^2$ :

\begin{theorem} \label{C1}
Let $u\in C(\rr  : L^2(\rr))$ be the global solution of the IVP \eqref{kgKdV} with $k=1$ provided by Theorem C
below. If there exist $\alpha>0$ and two different times $t_0,\,t_1\in\rr$ such that
\begin{equation}
\label{00.1}
|x|^{\alpha} u(x,t_0),\;|x|^{\alpha}u(x,t_1)\in L^2(\rr),
\end{equation}
then
\begin{equation}
\label{00.2}
u\in C(\rr  : H^{2\alpha}(\rr)).
\end{equation}
\end{theorem}

\begin{theorem} \label{D1}
Let $u,\,v\in C(\rr  : H^1(\rr))$ be global solutions of the IVP \eqref{kgKdV} with $k=1$ provided by Theorem C.
 If there exist $\alpha>1/2$ and two different times $t_0,\,t_1\in\rr$ such that
\begin{equation}
\label{000.1}
|x|^{\alpha} (u(x,t_0)-v(x,t_0)),\;|x|^{\alpha}(u(x,t_1)-v(x,t_1))\in L^2(\rr),
\end{equation}
then
\begin{equation}
\label{000.2}
u-v\in C(\rr  : H^{2\alpha}(\rr)).
\end{equation}
\end{theorem}

Combining Theorem \ref{A1} and Theorem \ref{C1} and taking an initial data $u_0\in L^2(\rr)$ with compact support such that 
$$
u_0\notin H^s(\rr),\;\;\;\forall \,s>0,
$$
one gets:

\begin{corollary} 
\label{corollary1}
There exists a solution 
\begin{equation}
\label{cc1}
u\in C(\rr : L^2(\rr)) \cap\dots . 
\end{equation}
of the IVP \eqref{kgKdV} with $k=1$ provided by Theorem C such that 
\begin{equation}
\label{cc2}
u(\cdot,0)=u_0(\cdot)
\end{equation}
has compact support and 
\begin{equation}
\label{cc3}
\begin{aligned}
&u(\cdot,t)\in C^{\infty}(\rr),\;\;\;\;\;\;\;\;\;\;\;\;\;\;\;\;\;\forall \,t \neq 0,\\
&u(\cdot,t)\notin L^2(|x|^{\epsilon}dx),\;\;\;\;\;\;\;\forall \,\epsilon>0,\;\;\forall \,t \neq 0,\\
&u(t)\in L^2(e^{a(t)x_{+}^{3/2}/2}dx),\;\;\;\;\;\;\;\,\forall \,t>0,\\
&u(t)\in L^2({e^{a(t)x_{-}^{3/2}/2}dx}),\;\;\;\;\;\;\;\,\forall \,t<0,
\end{aligned}
\end{equation}
\end{corollary}
with
\begin{equation}
 \label{cc4}
a(t)=\frac{a_0}{(1+27 a_0^2|t|/4)^{1/2}}.
 \end{equation}

 The results in Theorem \ref{C1} and Corollary \ref{corollary1} extend to the other powers $k=2,3,4...$ in \eqref{kgKdV} with the appropriate modifications,  accordingly to the value of $k$,  on  the regularity 
 and on  the length of the time interval $[0,T]$ as it was described in the statement of Theorem \ref{BB1}.
 
 The result in Theorem \ref{D1} holds for any power $k$ in \eqref{kgKdV} and any pair of solutions
 $$
 u,\;v\in C([-T,T]:H^1(\rr))
 $$
for any $T>0$ if $k=2,3$ and for $T=T(\|u_0\|_{1,2}, \|v_0\|_{1,2})$ for $k\geq 4$. 

 As a second consequence of our results above one has that there exist compact perturbations
 of the travelling wave solution with speed $c$
 \begin{equation}
 \label{solsol}
 u_{k,c}(x,t)=\phi_{k,c}(x-ct),
 \end{equation} 
 with 
 \begin{equation}
 \label{soliton}
 \phi_{k,c}(x)=(c_k \,c\,{\rm sech}^2(k\sqrt{c}\,x/2))^{1/k}
 \end{equation}
for the equation in \eqref{kgKdV} which destroy its exponential decay character:
 \begin{corollary} 
\label{corollary2}
For a given data of the form
\begin{equation}
\label{c5}
u_0(x)=\phi_{k,c}(x)+v_0(x),
\end{equation}
with $v_0\in H^1(\rr)$ compactly supported such that $v_0\notin H^{1+\epsilon}(\rr)$ for any $\epsilon>0$,
then the corresponding solution of the IVP provided by Theorem D
\begin{equation}
\label{c6}
u\in C([-T,T] : H^1(\rr)) \cap\dots
\end{equation}
of the IVP \eqref{kgKdV} satisfies that
\begin{equation}
\label{c7}
u(\cdot,t)\notin L^2(|x|^{1+\epsilon}dx),\;\;\;\;\;\;\;\forall \,\epsilon>0,\;\;\forall \,t \in [-T,T]-\{0\}.
\end{equation}
\end{corollary}

Remarks:

 (a) As in Corollary \ref{corollary1}, for $t>0$  the loss of decay is in the left hand side of $\rr$, and for $t<0$ in the right hand side of $\rr$. 
 
 (b) Combining the results in Theorem \ref{C1},  and its extension for all the equations in \eqref{kgKdV} commented above, with those found in \cite{NP12} one can also conclude that for $k=2,4,5,...$
 $$
 |x|^{\alpha}u(\cdot, t)\in L^2(\rr),\;\;\;\;\;\;\;\forall \,t\in[-T,T],
 $$
 and that for $k=1,3$ for any $\epsilon>0$
  $$
 |x|^{\alpha-\epsilon}u(\cdot, t)\in L^2(\rr),\;\;\,t\in[-T,T]-[t_0,t_1],\;\;\;\;\;|x|^{\alpha}u(\cdot, t)\in L^2(\rr),\;\;\;t\in[t_0,t_1].
 $$

 The main difference between the cases $k = 2, 4,5,...$ and $k = 1, 3$  is that for the later
the best available  well-posedness results require the use of the spaces $X_{s,b}$ defined in \eqref{xsb}, which makes fractional weights difficult to handle.
 
 (c) The equivalent of Theorem \ref{C1}  for the semi-linear Schr\"odinger equation \eqref{NLS} in all dimension $n$ was obtained in \cite{NP09}.

 We need to recall some results concerning the well-posedness (local and global) of the IVP \eqref{kgKdV}. First, we remember the definition of the space $X_{s,b}$ introduced in the context of dispersive equations in \cite{Bo}. 
  
  For $s, b\in\rr$, $X_{s,b}$ denotes the completion of the Schwartz space $\mathcal S(\rr^2)$ with respect to the norm
  \begin{equation}
  \label{xsb}
  \|F\|_{X_{s,b}}=(\int^{\infty}_{-\infty} \int^{\infty}_{-\infty} (1+|\tau-\xi^3|)^{2b}(1+|\xi|)^{2s}|\widehat{F}(\xi,\tau)|^2d\xi d\tau)^{1/2}.
  \end{equation}

The following result was established in \cite{Bo}, see also \cite{KPV96}:

   \begin{TC} \label{Bo}
 There exists $b\in(1/2,1)$ such that for any $u_0\in L^2(\rr)$ there exists a unique solution $u(t)$  of the IVP \eqref{kgKdV} with $k=1$ satisfying 
\begin{equation}
\label{1.09}
 u \in C([0,T]:L^2(\rr)),\;\;\;\;\;\text{for any}\;\;T>0,
 \end{equation} 
 \begin{equation}
 \label{1.10} 
 u\in X_{0,b},\;\;\;\partial_x(u^2)\in X_{0,b-1},\;\;\;\partial_t u\in X_{-3,b-1}.
 \end{equation}
Moreover, the map data-solution from  $L^2(\rr)$ into the class defined in \eqref{1.09}-\eqref{1.10} is Lipschitz for any $T>0$.
\end{TC}
 
 We recall that if $b>1/2$ one has, using Strichartz estimates and Kato local smoothing effects, that
 \begin{equation}
 \label{stri}
 \|f\|_{L^4([0,T]:L^{\infty}(\rr))}=(\int_0^T\|f(\cdot,t)\|^4_{\infty}dt)^{1/4}\leq c\|f\|_{X_{0,b}},
 \end{equation}
 and
  \begin{equation}
 \label{kat}
 \|\partial_x f\|_{L^{\infty}(\rr:L^2[0,T])}=
 \sup_{x\in\rr}(\int_0^T|\partial_xf(x,t)|^2dt)^{1/2}\leq c\|f\|_{X_{0,b}}.
 \end{equation}

 Therefore, combining these estimates and  Theorem C one has that the map data-solution is continuous from $L^2(\rr)$ to $ L^4([0,T]:L^{\infty}(\rr))$, for any $T>0$.
 
 Gathering the local and the global well-posedness results in \cite{KPV93}, \cite{CKSTT}, \cite{Ta}, \cite{GuPaSi} and  \cite{Gu} one has:
 
\begin{TD} \label{kgeq1} {\rm(a)} The IVP \eqref{kgKdV} with $k=2$ is globally well-posed in $H^s(\rr)$ for $s\geq 1/4$   (see \cite{CKSTT}).

{\rm(b)} The IVP \eqref{kgKdV} with $k=3$ is locally well-posed in $H^s(\rr)$ for $s\geq -1/6$, and globally  well-posed  in $H^s(\rr)$ for $s>-1/42$ (see \cite{Gu}, \cite{Ta}, \cite{GuPaSi}).

{\rm (c)} The IVP \eqref{kgKdV} with $k\geq 4$ is locally well-posed in $H^s(\rr)$ for $s\geq (k-4)/2k$ (see \cite{KPV93}).
\end{TD}
 
The rest of this paper is organized as follows:  In Section 2 Theorem \ref{A1} will be proved.  Since the proofs of Theorems \ref{B1}
and \ref{BB1}  are similar to the proof of Theorem \ref{A1} we will omit them.  The proof of Theorem \ref{C1} will be given  in Section 3.   
The proof of Theorem \ref{D1} will be omitted since its proof follows analogous arguments as the ones of Theorem \ref{C1}.  
Corollaries \ref{corollary1} and \ref{corollary2} are direct consequence of the previous results.

 \vskip.1in

 \section{Proof of Theorem \ref{A1}}\label{section2}

 Consider the IVP 
 \begin{equation}
\label{1.20}
\begin{cases}
\begin{aligned}
&\p_t u + \p_x^3 u + u\p_x u=0, \;\;\;\;t\geq 0,\;\;\;x\in  \rr,\\
&u(x,0)=u^{\epsilon}_0(x).\end{aligned}
\end{cases}
\end{equation}
 where
 \begin{equation}
 \label{1.21}
 u_0^{\epsilon}(x)=\rho_{\epsilon}\ast u_0(\cdot+\epsilon)(x)=\int_{-\infty}^{\infty}\,\frac{1}{\epsilon}\rho(\frac{y}{\epsilon})u_0(x+\epsilon-y)dy,
 \end{equation}
 where $\rho\in C^{\infty}(\rr)$  with $supp\,\rho\subset\{x\in\rr\,:\;|x|\leq 1\}$, $\,\rho\geq 0$ and $\,\int\rho(x)dx=1$. We claim that for any $\epsilon \in (0,1)$
 \begin{equation}
 \label{1.22}
 u_0^{\epsilon}\in H^{\infty}(\rr)\cap L^2(e^{a_0x_{+}^{3/2}}dx).
 \end{equation}
 To see \eqref{1.22} we use Minkowski's integral inequality to write:
 \begin{equation}
 \label{1.23}
 \begin{aligned}
& \|e^{a_0x_{+}^{3/2}/2}u_0^{\epsilon}\|_2=(\int e^{a_0x_{+}^{3/2}}|\int\frac{1}{\epsilon}\rho(\frac{y}{\epsilon})u_0(x+\epsilon-y)dy|^2dx)^{1/2}\\
 &=(\int |\int e^{a_0x_{+}^{3/2}/2}\rho_{\epsilon}(y)u_0(x+\epsilon-y)dy|^2dx)^{1/2}\\
 &\leq \int \rho_{\epsilon}(y)(\int e^{a_0x_{+}^{3/2}}|u_0(x+\epsilon-y)|^2dx)^{1/2}dy\\
 &\leq \int \rho_{\epsilon}(y)(\int e^{a_0(x-\epsilon+y)_{+}^{3/2}}|u_0(x)|^2dx)^{1/2}dy\\
& \leq \int \rho_{\epsilon}(y)(\int e^{a_0x_{+}^{3/2}}|u_0(x)|^2dx)^{1/2}dy=(\int e^{a_0x_{+}^{3/2}}|u_0(x)|^2dx)^{1/2},
 \end{aligned}
 \end{equation}
 since for  $y\in supp(\rho_{\epsilon})$ one has that $\,-\epsilon+y\leq 0$.  
 
 Also we recall  that 
 \begin{equation}
 \label{1.24}
 \lim_{\epsilon\downarrow 0}\|\,u_0^{\epsilon}-u_0\|_2=0.
 \end{equation}
 
 Therefore, for any $\epsilon\in (0,1)$ the corresponding solution $u^{\epsilon}(\cdot)$ of the IVP \eqref{1.20} has the properties stated  in Theorems A and C. 
 
 Next, for any $N\in \mathbb Z^{+}, \,N>>1$  we define the weight
 \begin{equation}
 \label{1.25}
 \varphi_N(x,t)=
 \begin{cases}
 \begin{aligned}
 &\;e^{a(t)/4},\;\;\;\;\;\;\;\,x\leq 0,\\
 &\;e^{a(t)\theta(x)},\;\;\;\;\;\,0\leq x\leq 1,\\
 &\;e^{a(t)x^{3/2}},\;\;\;\;\;1\leq x\leq N,\\
 &P_2(x,t),\;\;\;\;\;\;x\geq N,
 \end{aligned}
 \end{cases}
 \end{equation}
 where
 \begin{equation}
 \label{1.26}
 a'(t)+\frac{27}{8}a^{3}(t)=0,\;\;\;a(0)=a_0,
 \end{equation}
 i.e. 
 \begin{equation}
 \label{1.26b}
a(t)=\frac{a_0}{(1+27 a_0^2t/4)^{1/2}}\in (0,a_0],\;\;\;\;\;\;\forall \,t\geq 0,
 \end{equation}
 
 \begin{equation}
 \label{1.27}
 \theta(x)=\frac{1}{4}+\frac{15}{8}x^3-\frac{12}{8}x^4+\frac{3}{8}x^5,
 \end{equation}
and 
\begin{equation}
\label{1.28}
\begin{aligned}
P_2(x,t) &=e^{a(t)N^{3/2}}+\frac{3}{2}a(t)N^{1/2} e^{a(t)N^{3/2}}(x-N)\\
&+((\frac{3}{2}a(t)N^{1/2})^2+\frac{3}{4}a(t)N^{-1/2})e^{a(t)N^{3/2}}\frac{(x-N)^2}{2}.
\end{aligned}
\end{equation}

 Remarks :
 
 (i) For any $N\in \mathbb Z^{+}$, $x\geq 0$ and $t\geq 0$  one has
 \begin{equation}
 \label{new1}
 \varphi_N(x,t)\leq C(a_0)\, e^{a(t)x^{3/2}},
 \end{equation}
 
 (ii) the function $\theta(\cdot)$ matches the values of $x^{3/2} $ at $x=1$ and those ones of the constant function  $f(x)=1/4$ at $x=0$ and their derivatives up to order two
 and $P_2(x,t)$ matches those ones of $\,e^{a(t)x^{3/2}}$ at $\,x=N$ up to order two. Hence,  $\varphi_N(\cdot,t)\in C^2(\rr)$ and $\varphi_N(\cdot,t)\in C^3(\rr-\{0,1,N\})$
for all $t\geq 0$,

(iii) one has that
$$
 \;\;\;\;\; \theta''(x)=\frac{3x}{4}(15-24x+10x^2)=\frac{3x}{4}((\sqrt{10}x-\frac{12}{\sqrt{10}})^2+\frac{3}{5})\geq0,
$$
with $\theta'(0)=0$, hence $\theta'(x)\geq 0,\;\;\;x\in(0,1]$.

\vskip.1in

(iv) $\,P_2(x,t)\geq 0$  and
$$
\aligned
\partial_xP_2(x,t)&=
\frac{3}{2}a(t)N^{1/2}e^{a(t)N^{3/2}}\\
&+((\frac{3}{2}a(t)N^{1/2})^{2}
+\frac{3}{4}a(t)N^{-1/2})e^{a(t)N^{3/2}}(x-N)\geq 0,
\endaligned
$$
for $x\geq N$ and $t\geq 0$.

Therefore, 
\begin{equation}
\label{1.29}
\varphi_N(x,t),\;\;\partial_x\varphi_N(x,t)\geq 0,\;\;\;\;\;(x,t)\in\rr\times [0,\infty).
\end{equation}

Moreover, for $x\geq N$
\begin{equation}
\label{1.26bb}
\begin{aligned}
\partial_xP_2(x,t)&\leq \frac{3}{2}aN^{1/2}P_2(x,t)+\frac{3}{2}aN^{1/2}e^{aN^{3/2}}(x-N)(\frac{3}{2}aN^{1/2}+\frac{1}{2N})\\
&\leq  \frac{3}{2}aN^{1/2}P_2(x,t) +(\frac{3}{2}aN^{1/2}+\frac{1}{2N})P_2(x,t)\\
&\leq (1+3a(t)N^{1/2})P_2(x,t)\\
&\leq (1+3a_0x^{1/2})\varphi_N(x,t).
\end{aligned}
\end{equation}

Next, we multiply the equation in \eqref{1.20} by $u^{\epsilon} \varphi_N(x,t)$ where $u^{\epsilon}$ denotes the solution of the IVP \eqref{1.20}, integrating the result and 
formally using integration by parts one has
\begin{equation}
\label{3.1}
\begin{aligned}
\frac{d\;}{dt}\int (u^{\epsilon})^2\varphi_N dx&+3\int(\partial_x u^{\epsilon})^2\partial_x\varphi_Ndx\\
&=\int(u^{\epsilon})^2(\partial_x^3\varphi_N+\partial_t\varphi_N)dx+\frac{2}{3}\int(u^{\epsilon})^3\partial_x\varphi_Ndx.
\end{aligned}
\end{equation}

To justify the integration by parts used to obtain \eqref{3.1} we observe  that
by Theorems A and Theorem C
\begin{equation}
\label{3.2}
u^{\epsilon}\in C([0,\infty):H^{\infty}(\rr)),\;\;\;\;\;\;\;\;e^{\beta x}u^{\epsilon}\in C([0,\infty):L^{2}(\rr))\;\;\;\;\;\forall \,\beta>0.
\end{equation}

In general, if for some $\beta>0$ 
$e^{\beta x}f,\;\partial_x^2f\in L^2(\rr)$, then $\,e^{\beta x/2}\partial_xf \in L^2(\rr)$, and so $e^{\beta x/2}f \in H^1(\rr)$ (in particular, $e^{\beta x/2}f(x)\to 0\;\;\text{as}\;\;|x|\to \infty$) since 
\begin{equation}
\label{3.3}
\int e^{\beta x}(\partial_xf)^2dx\leq \beta^2\,\int e^{\beta x} f^2 dx +|\int e^{\beta x} f\,\partial_x^2fdx|.
\end{equation}
To prove \eqref{3.3} one first assumes that $f\in H^2(\rr)$ with compact support to obtain \eqref{3.3} by integration by parts, and then use the density of this class
to get the desired result. 

Thus, reapplying the last argument and using \eqref{3.2} it follows that for 
$$
e^{\beta x}\partial_x^ju^{\epsilon}\in L^2(\rr),\;\,j=1,2,3,
$$
and so for any $t\geq 0$ 
\begin{equation}
\label{3.5}
e^{\beta x}\partial^j_{x}u^{\epsilon}(x,t)\to 0\;\;\;j=0,1,2,\;\;\;\text{as}\;\;|x|\to \infty,
\end{equation}
which justify all integration by parts used to get \eqref{3.1} since at $+\infty$, $\varphi_N(x,t)$ as a function of $x$ is a polynomial of order two.

 Since 
 $$
 \partial_x\varphi_N(x,t)\geq 0,
 $$
 one can omit the second term on the left hand side (l.h.s.) of \eqref{3.1} to write
 \begin{equation}
\label{3.6}
\frac{d\;}{dt}\int (u^{\epsilon})^2\varphi_N dx\leq \int(u^{\epsilon})^2(\partial_x^3\varphi_N+\partial_t\varphi_N)dx+\frac{2}{3}\int(u^{\epsilon})^3\partial_x\varphi_Ndx.
\end{equation}

We shall consider the right hand side (r.h.s) of \eqref{3.6} in four different domains:
\begin{equation}
\label{3.7} 
\text{(a)}\;\;x\leq 0,\;\;\;\text{(b)}\;\;0\leq x\leq 1,\;\;\;\text{(c)}\;\;1\leq x\leq N,\;\;\;\text{(d)}\;\;x\geq N.
\end{equation}

  In the domain (a) ($x\leq 0$) one has
  $$
  \partial_x^j\varphi_N(x,t)=0\;\;\;j=1,2,3, \;\;\;\;\;\text{and}\;\;\;\;\partial_t\varphi_N(x,t)\leq 0,
  $$
  therefore the contribution of this domain to the r.h.s. of \eqref{3.6} is non-positive so it does not have to be considered.
  
  In the domain (b) ($0\leq x< 1$) one has 
  $$
  \aligned
  &\varphi_N(x,t)\leq \varphi_N(1,t)= e^{a(t)}\leq e^{a_0},\\
  &\partial_t\varphi_N(x,t)=a'(t)\varphi_N(x,t)\leq 0,\\
  &\partial_x\varphi_N(x,t)=a(t)\theta'(x)\varphi_N(x,t)\leq c \,a_0\varphi_N(x,t),\\
  &\partial_x^2\varphi_N(x,t)=(a(t)\theta^{''}(x)+(a(t)\theta'(x))^2)\varphi_N(x,t)\leq c\,(a_0+a_0^2)\varphi_N(x,t),\\
  &\partial_x^3\varphi_N(x,t)=(a(t)\theta^{(3)}(x)+3(a(t))^2\theta^{''}(x)\theta'(x)+(a(t)\theta'(x))^{(3)})\varphi_N(x,t)\\
  &\;\;\;\;\;\;\;\;\;\;\;\;\;\;\;\;\;\leq c\,(\,a_0 +a_0^2+a_0^3)\varphi_N(x,t).
\endaligned
$$  
Therefore
\begin{equation}
\label{3.8}
\int_0^1(u^{\epsilon})^2(\partial_x^3\varphi_N+\partial_t\varphi_N)dx \leq K(a_0)\int_0^1(u^{\epsilon})^2\varphi_N dx,
\end{equation}
and
\begin{equation}
\label{3.9}
\int_0^1(u^{\epsilon})^3\partial_x\varphi_N dx \leq K(a_0)\,\|u^{\epsilon}(t)\|_{L^{\infty}(0\leq x\leq 1)}\int_0^1(u^{\epsilon})^2\varphi_N dx.
\end{equation}

  In the domain (c) ($1\leq x< N$) one has 
  \begin{equation}
  \label{3.10}
  \partial_x^3\varphi_N+\partial_t\varphi_N=(\frac{27}{8}a^3x^{3/2}+\frac{27}{8}a^2-\frac{3}{8}ax^{-3/2}+a'x^{3/2})\varphi_N(x,t).
  \end{equation}
For the third term on the r.h.s. of \eqref{3.10}  we observe that
\begin{equation}
\label{3.10b}
 -\frac{3}{8}a(t)x^{-3/2}e^{a(t)x^{3/2}}\leq 0,
 \end{equation}
 so it can be omitted. Next, we use that
 \begin{equation}
 \label{3.11}
 a'(t)+\frac{27}{8}a^3(t)=0,
 \end{equation}
 which combined with \eqref{3.10b} and \eqref{3.11} gives us the estimate
 \begin{equation}
 \label{3.12}
 \int_1^N(u^{\epsilon})^2(\partial_x^3\varphi_N+\partial_t\varphi_N)dx\leq \frac{27}{8}a^2(t)\int_1^N(u^{\epsilon})^2\varphi_Ndx
 \leq \frac{27}{8}a_0^2 \int_1^N(u^{\epsilon})^2\varphi_Ndx.
 \end{equation}

For the other term on the r.h.s. of \eqref{3.6} in this domain one sees that
\begin{equation}
\label{3.13}
\begin{aligned}
\int_1^N(u^{\epsilon})^3\partial_x\varphi_Ndx&=\frac{3}{2}\int_1^Na(t)x^{1/2}(u^{\epsilon})^3\varphi_Ndx\\
&\leq \frac{3}{2}a_0\|x^{1/2}u^{\epsilon}(t)\|_{L^{\infty}(1<x<N)}\int_1^N(u^{\epsilon})^2\varphi_Ndx.
\end{aligned}
 \end{equation}

 Finally, we consider the domain (d) ($x\geq N$).  Since in this domain 
 $$
 \partial_x^3P_2(x,t)\equiv 0,\;\;\;\;\;\;\partial_tP_2(x,t)=a'(t)\{\cdot\}\leq 0,
 $$
one has  
 $$
  \int_N^{\infty}(u^{\epsilon})^2(\partial_x^3\varphi_N+\partial_t\varphi_N)dx\leq 0,
  $$
 so the contribution of this term on the r.h.s. of \eqref{3.6} does not need to be considered. To bound the contribution of the  other term on the r.h.s. of \eqref{3.6} 
  using \eqref{1.26bb} we write 
\begin{equation}
\label{3.14}
\begin{aligned}
\int_N^{\infty}(u^{\epsilon})^3&\partial_x\varphi_Ndx\\
&=\int_N^{\infty}(u^{\epsilon})^3(1+3a_0x^{1/2})P_2(x,t)dx\\
&\leq K(a_0)\,\|x^{1/2}u^{\epsilon}(t)\|_{L^{\infty}(x\geq N)} \,\int_N^{\infty}(u^{\epsilon})^2\varphi_N(x,t)dx.
\end{aligned}
\end{equation}

  It remains to estimate the terms in the $L^{\infty}$-norm in \eqref{3.9},  \eqref{3.13} and \eqref{3.14}. We estimate the terms \eqref{3.13} and \eqref{3.14}, the estimation in \eqref{3.9} being similar.    For  \eqref{3.13} and \eqref{3.14}, using   \eqref{1.6} in Theorem A and an estimate similar to that in \eqref{1.23}, we have 
  \begin{equation}
  \label{3.16}
  \begin{aligned}
  \|x^{1/2}u^{\epsilon}(t)\|_{L^{\infty}(x>1)}&\leq \|e^x u^{\epsilon}(t)\|_{\infty}\leq  \|e^x u^{\epsilon}(t)\|^{1/2}_2\| \partial_x(e^x u^{\epsilon}(t))\|^{1/2}_2\\
  &\leq  \|e^x u^{\epsilon}(t)\|^{1/2}_2( \|e^x u^{\epsilon}(t)\|^{1/2}_2+\|e^x \partial_xu^{\epsilon}(t)\|^{1/2}_2)\\
  &\leq  2\|e^x u^{\epsilon}(t)\|_2+\|e^x \partial_xu^{\epsilon}(t)\|_2\\
  &\leq 2e^M\|e^xu_0^{\epsilon}\|_2 + \|e^x \partial_xu^{\epsilon}(t)\|_2\\
  &\leq 2e^M\|e^xu_0\|_2 + \|e^x \partial_xu^{\epsilon}(t)\|_2.
  \end{aligned}
  \end{equation}
 
 From \eqref{1.7} in Theorem A
 and an argument similar to that in \eqref{1.23} we obtain the bound for the integral in time interval $[0,T]$ of the last term in \eqref{3.16} 
 \begin{equation}
 \label{3.17}
 \begin{aligned}
 \int_0^T\| e^x\partial_xu^{\epsilon}(t)\|_2dt &\leq T^{1/2}(\int_0^T\| e^x\partial_xu^{\epsilon}(t)\|^2_2dt)^{1/2}\\
& \leq T^{1/2}\frac{1}{4} e^{MT}\| e^x u^{\epsilon}_0\|_2^2  \leq T^{1/2}\frac{1}{4} e^{MT}\| e^x u_0\|_2^2,
 \end{aligned}
  \end{equation}
   with $M=M(\|u_0\|_2)$.

  Inserting the above estimates in \eqref{3.6}, using Gronwall's inequality, and applying \eqref{new1} with $t=0$,  one gets that
  \begin{equation}
  \label{3.18}
\sup_{t\in[0,T]} \int_{-\infty}^{\infty}|u^{\epsilon}(x,t)|^2\varphi_N(x,t)dx\leq C^*=C^*( a_0, \|u_0\|_2, \|e^{a_0x_{+}^{3/2}/2}u_0\|_2, T).
\end{equation}
  
  We observe that for each $N\in\mathbb Z^{+}$ fixed there exists $c_N$ such  that
  \begin{equation}
  \label{3.19}
  \varphi_N(x,t)\leq c_Ne^{a_0x},\;\;\;\;\;\;\;\,x\geq 0,
  \end{equation}
and by the continuity of the map data-solution, in Theorem A   it follows that
\begin{equation}
\label{3.20}
\sup_{t\in[0,T]}\int_{-\infty}^{\infty}e^{a_0x}|u^{\epsilon}(x,t)-u(x,t)|^2dx\to 0\;\;\;\;\text{as}\;\;\;\;\epsilon\downarrow 0.
\end{equation}
Combining \eqref{3.19} and \eqref{3.20} one concludes that
\begin{equation}
\label{3.21}
\sup_{t\in[0,T]}\int_{-\infty}^{\infty}\varphi_N(x,t)|u^{\epsilon}(x,t)-u(x,t)|^2dx\to 0\;\;\;\;\text{as}\;\;\;\;\epsilon\downarrow 0.
\end{equation}
Hence, for any $t\in[0,T]$ 
$$
\int_{-\infty}^{\infty}\varphi_N(x,t)|u^{\epsilon}(x,t)|^2dx\to \int_{-\infty}^{\infty}\varphi_N(x,t)|u(x,t)|^2dx\;\;\;\;\text{as}\;\;\;\;\epsilon\downarrow 0,
$$
and consequently, for any $t\in [0,T]$
\begin{equation}
\label{3.22}
\int_{-\infty}^{\infty}\varphi_N(x,t)|u(x,t)|^2dx \leq C^*=C^*( a_0, \|u_0\|_2, \|e^{a_0x_{+}^{3/2}/2}u_0\|_2, T).
\end{equation}

Finally, letting $N\uparrow \infty$ in \eqref{3.22}  Fatou's lemma yields the desired result
\begin{equation}
\label{3.22b}
\int_{-\infty}^{\infty} e^{a(t)x_{+}^{3/2}}|u(x,t)|^2dx \leq C^*=C^*( a_0, \|u_0\|_2, \|e^{a_0x_{+}^{3/2}/2}u_0\|_2, T).
\end{equation}

\vskip.2in

 \section{Proof of Theorem \ref{C1}}\label{section3}
 
 Without loss of generality we assume $t_0=0$ and $t_1>0$.
 
First we shall consider the case $\alpha\in (0,1/2]$.
 
 For $x\geq 0$, $N\in\mathbb Z^+$ and $\alpha>0$ we define
 \begin{equation}
 \label{5.1}
 \varphi_{N,\alpha}(x)=
 \begin{cases}
 \begin{aligned}
 &(1+x^4)^{\alpha/2}-1,\;\;\;\;\;\;\;x\in[0,N],\\
 &(2N)^{2\alpha},\,\;\;\;\;\;\;\;\;\;\;\;\;\;\;\;\;\;\;x\geq 10N,
 \end{aligned}
 \end{cases}
 \end{equation}
 with $\varphi_{N,\alpha}\in C^3((0,\infty))$, $\varphi_{N,\alpha}(x)\geq 0$, $\varphi^{'}_{N,\alpha}\geq 0$ and
 for $\alpha\in (0,1/2]$
  \begin{equation}
 \label{5.2}
 |\varphi_{N,\alpha}^{(j)}(x)|\leq c,\;\;\;j=1,2,3, \;\;\;\;\;c\;\text{ independent of }\;N.
 \end{equation}
 Let $\,\phi_{N,\alpha}\,$ be the odd extension of $\,\varphi_{N,\alpha}$, i.e. 
  \begin{equation}
 \label{5.3}
 \phi_{N,\alpha} (x)=\phi_N(x)=
 \begin{cases}
 \begin{aligned}
 &\;\;\;\;\varphi_{N,\alpha}(x),\;\;\;\;\;\;\;\;\;\;x\geq 0,\\
 &-\varphi_{N,\alpha}(-x),\;\;\;\;\;\;\;x\leq 0.
 \end{aligned}
 \end{cases}
 \end{equation}
 Notice that
 \begin{equation}
 \label{5.4}
 \phi_N^{'}(x)\geq 0,\;\;\;\forall \,x\in\rr,\;\;\;\phi_N\in C^3(\rr),\,\;\;\, \|  \phi_N\|_{\infty}=(2N)^{2\alpha}.
 \end{equation}
 
 Next, we consider a sequence of data $(u_{0,m})_{m\in\mathbb Z^{+}}\subset  \mathcal S(\rr)$ such that
 \begin {equation}
 \label{5.5}
 \|u_0-u_{0,m}\|_2\to 0\;\;\;\;\;\text{as}\;\;\;\;m\uparrow \infty,
 \end{equation} 
 and denote by $u_m$ the solution of the IVP \eqref{kgKdV} with $k=1$ and data $u_m(x,0)=u_{0,m}(x)$. From the results in \cite{Ka} 
 one has that
 \begin{equation}
 \label{5.6}
 u_m\in C(\rr : \mathcal S(\rr)).
 \end{equation}
 By the continuous dependence
 of the solution upon the data (see Theorem C and comments after its statement) one has that for each $T>0$ 
 \begin{equation}
 \label{5.7}
 \begin{aligned}
 (a)\;&\sup_{t\in[0,T]}\|u(t)-u_{m}(t)\|_2\to 0\;\;\;\;\;\text{as}\;\;\;\;m\uparrow \infty,\\
 (b)\;&\int_{-T}^T\|u(t)-u_{m}(t)\|^4_{\infty}dt \to 0\;\;\;\;\;\text{as}\;\;\;\;m\uparrow \infty,\\
 (c)\;&\sup_{x\in\rr}\,\int_{-T}^T |\partial_x(u-u_m)(x,t)|^2dt\to 0\;\;\;\;\;\text{as}\;\;\;\;m\uparrow \infty.
 \end{aligned}
 \end{equation}
 
 Since $u_m\in C(\rr : \mathcal S(\rr))$ satisfies the equation in \eqref{kgKdV} with $k=1$, multiplying it by $u_m\,\phi_N$ after integration by parts (justified since $\phi_N$ is bounded)
 one gets 
 \begin{equation}
 \label{5.8}
 \frac{d\;}{dt}\int (u_m)^2\phi_Ndx + 3\int (\partial_xu_m)^2\phi_N^{'}dx=\int (u_m)^2\phi_N^{(3)} dx  +\frac{2}{3}\int (u_m)^3\phi_N^{'}dx.
 \end{equation}
 We observe that for $m$ large enough 
 \begin{equation}
 \label{5.9}
 \begin{aligned}
 &|\int (u_m)^2\phi_N^{(3)} dx|\leq c\|u_{0,m}\|_2^2\leq 2c\|u_{0}\|_2^2,\\
 &|\int (u_m)^3\phi_N^{'}dx|\leq c\|u_m(t)\|_{\infty}\|u_{0,m}\|_2^2\leq 2c\|u_m(t)\|_{\infty}\|u_{0}\|_2^2.
 \end{aligned}
 \end{equation}
 
 Integrating in the time interval $[0,t_1]$ the identity \eqref{5.8} and using \eqref{5.9} it follows that
 \begin{equation}
 \label{5.10}
 \begin{aligned}
 \int_0^{t_1}\int (\partial_xu_m)^2(x,t)\phi_N^{'}(x)dxdt&\leq \|(u_m(t_1))^2\phi_N\|_1+\|(u_{0,m})^2\phi_N\|_1\\
 &+c\,t_1\|u_0\|_2^2+c\,t_1^{3/4}\,\|u_0\|^2_2\,(\int_0^{t_1}\|u_m(t)\|^4_{\infty}dt)^{1/4},
 \end{aligned}
 \end{equation}
 where $c$ denotes a constant whose value may change from line to line and is independent of the initial parameters of the problem. Letting $m\uparrow \infty$, 
 using \eqref{5.7} part (c), Theorem C, and \eqref{stri}, one gets that
 \begin{equation}
 \label{5.11}
 \begin{aligned}
\varlimsup_{m\uparrow \infty} \int_0^{t_1}\int (\partial_xu_m)^2(x,t)\phi_N^{'}&(x)dxdt\leq \|(u(t_1))^2\phi_N\|_1+\|(u_{0})^2\phi_N\|_1\\
 +&c\,t_1\|u_0\|_2^2+c\,t_1^{3/4}\,\|u_0\|^2_2\,(\int_0^{t_1}\|u(t)\|^4_{\infty}dt)^{1/4} \leq M,
 \end{aligned}
 \end{equation} 
with $M=M(\|\langle x\rangle^{\alpha} u_0\|_2,\,\|\langle x\rangle^{\alpha} u(t_1)\|_2)$.
Next, we use \eqref{5.7} part (c) to conclude that for any $\bar N\in \mathbb Z^{+}$ fixed
\begin{equation}
\label{5.12}
\partial_xu_m\to \partial_xu\;\;\;\;\;\text{in}\;\;\;L^2([-T,T]\times [-\bar N,\bar N] )\;\;\;\;\text{as}\;\;\;\;m\uparrow \;\infty.
\end{equation}

 Since $\phi^{'}_N$ has compact support one gets that
 \begin{equation}
 \label{5.13}
 \int_0^{t_1}\int (\partial_xu)^2(x,t)\phi_N^{'}(x)dxdt\leq M.
 \end{equation}
 Finally, we recall that  $\phi_N^{'}(x)$ is even, $\phi_N^{'}(x)\geq 0$, and for $x>1$
 $$
 \phi_N^{'}(x)\to \frac{2\alpha x^3}{(1+x^4)^{1-\alpha/2}}\sim \langle x\rangle^{2\alpha-1},
 $$
therefore using Fatou's lemma in \eqref{5.13} one concludes that
 \begin{equation}
 \label{5.14}
 \int_0^{t_1}\int_{|x|\geq 1} (\partial_xu)^2(x,t)\langle x\rangle^{2\alpha-1}\,dxdt\leq M.
 \end{equation}
 Since (see Theorem C and \eqref{kat})
 $$
  \int_0^{t_1}\int_{|x|\leq 1} (\partial_xu)^2(x,t)\,dxdt\leq M,
  $$
  one concludes that
  \begin{equation}
 \label{5.15}
 \int_0^{t_1}\int(\partial_xu)^2(x,t)\,\langle x\rangle^{2\alpha-1}\,dxdt\leq M.
 \end{equation}
 
 Once we have obtained \eqref{5.15} we reapply the above argument with $\psi_{N,\alpha}(x)=\psi_N(x)$ the even extension of $\varphi_{N,\alpha}$ instead of $\phi_N(x)$, i.e. 
 \begin{equation}
 \label{5.3b}
 \psi_{N,\alpha} (x)=\psi_N(x)=
 \begin{cases}
 \begin{aligned}
 &\;\;\;\;\varphi_{N,\alpha}(x),\;\;\;\;\;\;\;\;\;\;x\geq 0,\\
 &\;\varphi_{N,\alpha}(-x),\;\;\;\;\;\;\;\;\;\;x\leq 0.
 \end{aligned}
 \end{cases}
 \end{equation}
 We observe that 
 \begin{equation}
 \label{5.4b}
 |\psi_N'(x)|\leq c\langle x\rangle^{2\alpha-1}.
 \end{equation}
 
 Using the formula \eqref{5.8}  with $\psi_N(x)$ instead of $\phi_N(x)$, and estimates similar to that in \eqref{5.9} it follows that 
 $$
 \lim_{m\uparrow \infty} \int_0^{t_1}\int(\partial_xu_m)^2\psi'_N(x)dxdt=\int_0^{t_1}\int (\partial_xu)^2\psi'_N(x)dxdt,
 $$
 and
 $$
 |\int_0^{t_1}\int (\partial_xu)^2\psi'_N(x)dxdt|\leq \int_0^{t_1}\int (\partial_xu)^2\langle x\rangle^{2\alpha-1} dxdt\leq M.
 $$
 
 From  these estimates and integrating in the time interval 
   $[0,t]\subset [0,t_1]$ one  obtains that
  \begin{equation}
 \label{5.16}
 \langle x\rangle^{\alpha}\,u(t)\in L^2(\rr)\;\;\;\;\;\;\;t\in[0,t_1].
 \end{equation} 
 
Combining \eqref{5.15} and \eqref{5.16} one sees that there exists $t^*\in [0,t_1]$ such that
 $$
 \langle x\rangle^{\alpha}\,u(t^*),\;\partial_xu(t^*)\langle x\rangle^{\alpha-1/2}\in L^2(\rr).
 $$
 Hence, 
 \begin{equation}
 \label{5.17}
  \langle x\rangle^{\alpha}u(t^*),\;\;J(\langle x\rangle^{\alpha-1/2}u(\cdot,t^*))=(1-\partial_x^2)^{1/2}(\langle x\rangle^{\alpha-1/2}u(t^*))\in L^2(\rr).
  \end{equation}
 
  Next we shall use the following lemma (see \cite{BeLo} and  \cite{NP09}) :
 \begin{lemma}
\label{lemma1}
Let $a,\,b>0$. Assume that $ J^af=(1-\partial^2)^{a/2}f\in L^2(\rr)$ and \newline $\langle x\rangle^{b}f=
(1+|x|^2)^{b/2}f\in L^2(\rr)$. Then for any $\theta \in (0,1)$
\begin{equation}
\label{complex}
\|J^{ \theta a}(\langle x\rangle^{(1-\theta) b} f)\|_2\leq c\|\langle x\rangle^b f\|_2^{1-\theta}\,\|J^af\|_2^{\theta}.
\end{equation}
\end{lemma}
 
 Defining 
 $$
 f=\langle x\rangle^{\alpha-1/2}u(t^*)
 $$
 it follows that
 $$
 Jf,\;\langle x\rangle^{1/2} f\in L^2(\rr)
 $$
 and one gets from \eqref{complex} with $\theta=2\alpha$ that
 \begin{equation}
 \label{5.18}
 J^{2\alpha}(\langle x\rangle^{1/2-\alpha}f)=J^{2 \alpha}u(\cdot,t^*)\in L^2(\rr).
 \end{equation}
  Once \eqref{5.18} has been established the rest of the proof of Theorem \ref{C1} follows the argument described in \cite{NP12}.

  Next, we consider the case $\alpha\in (1/2,1]$.

 From the previous case we already know that
 $$
 u\in C([0,t_1] : H^1(\rr)\cap L^2(\langle x\rangle dx)).
 $$
 We also observe that for $\alpha\in (1/2,1]$
$$ 
  \phi_{N,\alpha}'(x) +|\psi^{'}_{N,\alpha}(x)| +|\phi_{N,\alpha-1/2}(x)|\leq c\,\langle x\rangle,\;\;\;\;\;\;\;\;\;| \phi_{N,\alpha}^{(3)}(x)|\leq c.
  $$
  
  As before we get 
   \begin{equation}
 \label{5.19}
 \frac{d\;}{dt}\int (u_m)^2\phi_{N,\alpha}dx + 3\int (\partial_xu_m)^2\phi_{N,\alpha}^{'}dx=\int (u_m)^2\phi_{N,\alpha}^{(3)} dx  +\frac{2}{3}\int (u_m)^3\phi_{N,\alpha}^{'}dx,
 \end{equation}
with $\, u_m\in C(\rr : \mathcal S(\rr))$ which justifies the integration by parts in \eqref{5.19}. 

In \eqref{5.19} we first use that
\begin{equation}
\label{5.20}
|\int (u_m)^2 \phi_{N,\alpha}^{(3)} dx|\leq c\|u_{0,m}\|^2_2\leq 2c\|u_0\|^2_2.
\end{equation}

 Next, we estimate the cubic  term in \eqref{5.19} (last term there).  For this we observe that
 \begin{equation}
 \label{5.21}
 |\phi_{N,\alpha}^{'}(x)| \leq c\, |\phi_{N,\alpha-1/2}(x)|^2,
 \end{equation}
  with $c$ independent of $N$. Thus,
  \begin{equation}
 \label{3.33}
  |\int  (u_m)^3 \phi_{N,\alpha}'(x)dx| \leq c \|u_m(t)\|_{\infty}\| u_m(t)\phi_{N,\alpha-1/2}\|_2^2.
  \end{equation}
  
  Combining the facts that for each $N$  fixed the $\phi_{N,\alpha}$'s are  bounded and that
  $$
  \sup_{[0,t_1]}\|(u_m-u)(t)\|_2\to 0\;\;\;\;\;\;\text{as}\;\;\;\,m\uparrow \infty,
  $$   
  one has 
 $$
  \sup_{[0,t_1]}\|(u_m-u)(t)\phi_{N,\alpha-1/2}\|_2\to 0\;\;\;\;\;\;\text{as}\;\;\;\,m\uparrow \infty,
  $$   
  and consequently,
  \begin{equation}
  \label{5.23}
   \sup_{[0,t_1]}\|u_m(t)\phi_{N,\alpha-1/2}\|_2\leq 
  2 \sup_{[0,t_1]}\|u(t)\phi_{N,\alpha-1/2}\|_2
  \leq  2\sup_{[0,t_1]}\|\langle x\rangle^{1/2} u(t)\|_2 \leq M,
  \end{equation}
  with $M =M(\|\langle x\rangle^{1/2}u_0\|_2,\,\|\langle x\rangle^{1/2}u(t_1)\|_2)$ for $m>>1$.
  
  Inserting the above estimates in \eqref{5.19} and following the argument in the previous case $\alpha\in(0,1/2]$ one gets that
 \begin{equation}
 \label{5.24}
 \int_0^{t_1}\int (\partial_x u_m)^2(x,t)\phi^{'}_{N}(x)dx dt\leq (1+t_1)M,
 \end{equation}
 for $m>>1$  and 
\begin{equation}
 \label{5.25}
 \int_0^{t_1}\int (\partial_x u)^2(x,t)\,\langle x\rangle^{2\alpha -1}dx dt\leq (1+t_1)M,
 \end{equation}  
  where $M =M(\|\langle x\rangle^{1/2}u_0\|_2,\,\| \langle x\rangle^{1/2}u(t_1)\|_2)$.
   
  Next, we deduce \eqref{5.19} with $\psi_{N,\alpha}$ instead of $
  \varphi_{N,\alpha}$ and use \eqref{5.24} and \eqref{5.25} as in the previous case to get that
  \begin{equation}
  \label{5.26}
  \langle x\rangle^{\alpha}u(t)\in L^2(\rr)\;\;\;\;\;\;t\in [0,t_1].
  \end{equation}
  
  Once that the estimates \eqref{5.25} and \eqref{5.26} have been established the argument in the previous case involving Lemma \ref{lemma1} provides the desired result.

   A similar  boot-strap argument can be used in the case $\alpha\in (1,3/2]$ and respectively for higher $\alpha$.

\vspace{3mm} \noindent{\large {\bf Acknowledgments}}
\vspace{3mm}\\  P. I. was supported by DIME Universidad Na\-cio\-nal de Co\-lom\-bia-Me\-de\-ll\'in, grant 
2010\-1001\-1032. F. L. was partially supported
by CNPq and FAPERJ/Brazil. G. P. was  supported by a NSF grant  DMS-1101499.

\end{document}